\documentclass[11pt]{article}\UseRawInputEncoding
\usepackage{amssymb,amsfonts,amsmath,latexsym,epsf,tikz,url}
\newtheorem{theorem}{Theorem}[section]

\newtheorem{corollary}[theorem]{Corollary}

\newcommand{\proof}{\noindent{\bf Proof.\ }}
\newcommand{\qed}{\hfill $\square$\medskip}

\begin{document}
\title{More on total domination polynomial and $\mathcal{D}_t$-equivalence classes of some graphs }

\author{Saeid Alikhani$^{}$\footnote{Corresponding author}  \and  Nasrin Jafari }

\date{\today}

\maketitle

\begin{center}

   Department of Mathematics, Yazd University, 89195-741, Yazd, Iran\\
{\tt alikhani@yazd.ac.ir, nasrin7190@yahoo.com}\\

\end{center}

\begin{abstract}
Let $G = (V, E)$ be a simple graph of order $n$. The total dominating set of $G$ is a subset $D$ of $V$ that every vertex of $V$ is adjacent to some vertices of $D$. The total domination number of $G$ is equal to minimum cardinality of  total dominating set in $G$ and is denoted by $\gamma_t(G)$. The total domination polynomial of $G$ is the polynomial $D_t(G,x)=\sum_{i=\gamma_t(G)}^n d_t(G,i)x^i$, where $d_t(G,i)$ is the number of total dominating sets of $G$ of size $i$. Two graphs $G$ and $H$ are said to be total dominating equivalent or simply $\mathcal{D}_t$-equivalent, if $D_t(G,x)=D_t(H,x)$. The equivalence class of $G$, denoted $[G]$, is the set of all graphs $\mathcal{D}_t$-equivalent to $G$.   In this paper, we investigate $\mathcal{D}_t$-equivalence classes of some graphs. Also we introduce some families of graphs whose total domination polynomials are unimodal.
\end{abstract}


\medskip
\noindent{\bf AMS Subj.\ Class.}: 05C30, 05C69.  

\section{Introduction}

Let $G = (V, E)$ be a simple graph. The order of $G$ is the number of vertices of $G$. For any vertex $ v \in V$, the open neighborhood of $v$ is the set $N(v)=\{ u \in V | uv \in E\}$ and the closed neighborhood is the set $N[v]=N (v) \cup \{v\}$.
For a set $S\subset V$, the open neighborhood of $S$ is the set $N(S)=\bigcup_{v\in S }N(v)$ and the closed neighborhood of $S$ is the set $N[S]=N (S) \cup S$. The set $D\subset V$ is a total dominating set if every vertex of $V$ is adjacent to some vertices of $D$, or equivalently, $N(D)=V$. The total domination  number $\gamma_t(G)$ is the minimum cardinality of a total dominating set in $G$. A total dominating set with cardinality $\gamma_t(G)$ is called a $\gamma_t$-set. An $i$-subset of $V$ is a subset of $V$ of cardinality $i$. Let $D_t(G, i)$ be the family of total dominating sets of $G$ which are $i$-subsets and let $d_t(G,i)=|D_t(G, i)|$. The polynomial $D_t(G, x)=\sum_{i=1}^n d_t(G,i)x^i$ is defined as total domination polynomial of $G$. A root of $D_t(G, x)$ is called a total  domination root of $G$.
For many graph polynomials, their roots have attracted considerable attention. 

A natural question to ask is to what extent can a graph polynomial describe the underlying graph
(for example, a survey of what is known with regards to chromatic polynomials can be found in Chapter 3 of \cite{Dong}). We say that two graphs $G$ and $H$ are total domination equivalent or simply $\mathcal{D}_t$-equivalent (written
$G\sim_{t}H$) if they have the same total domination polynomial. Similar to domination polynomial \cite{euro,clique}, we let $[G]$ denote the $\mathcal{D}_t$-equivalence class determined by $G$, that is $[G]=\{H|H\sim_{t} G\}$. A graph $G$ is said to be total dominating unique or simply $\mathcal{D}_t$-unique if $[G]=\{G\}$.
Two problems arise:
\begin{enumerate}
	\item [(i)] 
 Which graphs are $\mathcal{D}_t$-unique, that is, are completely determined by their total domination polynomials? 
\item[(ii)] Determine the $\mathcal{D}_t$-equivalence class of a graph? 
\end{enumerate}   
Both problems appear difficult, but there are some partial results known. 

Recurrence relations of graph polynomials have received considerable attention in the literature. It is well-known that the independence
polynomial and matching polynomial of a graph satisfies a linear recurrence relation with respect to two vertex elimination operations, the deletion of a vertex and the deletion of vertex's closed neighborhood. Other graph polynomials in the literature satisfy similar recurrence relations with respect to vertex and edge elimination operations \cite{Kotek}. 
In contrast, it is significantly harder to find recurrence relations for the domination polynomial and the total domination polynomial.  The easiest recurrence relation is to remove an edge and to compute the total domination polynomial of the graph arising instead of the one for the original graph. Indeed, for the total domination polynomial of a graph there might be such irrelevant edges, that can be deleted without changing the value of the total domination polynomial at all.  
An irrelevant edge is an edge $e \in E$ of $G$, such that $D_t(G, x) = D_t(G\setminus e, x)$. These edges can be useful to classify some graphs by their total domination polynomials. 

A finite sequence of real numbers $(a_0, a_1, a_2, \ldots , a_n)$ is said to be unimodal, if there is some $k \in \{0, 1, \ldots , n\}$, called the mode of sequence, such that

$$a_0 \leq \ldots \leq a_{k-1} \leq a_k\geq a_{k+1} \geq \ldots \geq a_n.$$

The mode is unique if $a_{k-1} < a_k > a_{k+1}$. A polynomial is called unimodal if the sequence of its coefficients is unimodal.
It is log-concave if $a_k^2=a_{k-1}a_{k+1}$ for all $1 \leq k \leq n - 1$. It is symmetric if $a_k = a_{n-k}$ for $0 \leq k\leq n$. A log-concave sequence of positive numbers is unimodal (see, e.g., \cite{Brenti1,Brenti2, Levit2}). We say that a polynomial is unimodal (log-concave, symmetric, respectively) if the sequence of its coefficients is unimodal (log-concave, symmetric, respectively). A mode of
the sequence $a_0, a_1, \ldots, a_n$ is also called a mode of the polynomial $\sum_{k=0}^n a_kx^k$.

Unimodality problems of graph polynomials have always been of great interest to researchers in graph theory \cite{Fjafari, Levit1}. There are a number of results concerning the coefficients of independence polynomials, many of which consider graphs formed by applying some sort of operation to simpler graphs. In \cite{Rosenfeld}, for instance, Rosenfeld examines the independence polynomials of graphs formed by taking various rooted products of simpler graphs(See \cite{Godsil} for the definition of the rooted product of two graphs.) In particular, he shows that the property of having all real roots is preserved under forming rooted products. Mandrescu in \cite{Mandrescu} has shown that the independence polynomial of corona product of any graph with $2$ copies of $K_1$, i.e., $I(G\circ 2K_1, x)$ is unimodal. Levit and Mandrescu in \cite{Levit2} generalized this result and have shown that if $H = K_r -e$, $r\geq 2$, then the polynomial $I(G \circ H, x)$ is unimodal and symmetric for every graph $G$.

Although the unimodality of independence polynomial has been actively studied, almost no attention has been given to the unimodality of total domination polynomials. We checked the total domination polynomial of graphs of order at most six (see \cite{atlas}) and observed that all of these polynomials are unimodal.  

As usual we denote the complete graph, path and cycle of order $n $ by $K_n$, $P_n$ and $C_n$, respectively. Also $S_n$ is the star graph with $n$ vertices.

The corona of two graphs $G_1$ and $G_2$, as defined by Frucht and Harary in \cite{harary}, is the graph $G= G_1\circ G_2$ formed from one copy of $G_1$ and $|V(G_1)|$ copies of $G_2$, where the i-th vertex of $G_1$ is adjacent to every vertex in the $i$-th copy of $G_2$. The corona $G\circ K_1$, in particular, is the graph constructed from a copy of $G$, where for each vertex $v \in V(G)$, a new vertex $v'$ and a pendant edge $vv'$ are added. 

\medskip
In the next section, we study the total dominating equivalence classes of some graphs such as $G\circ \overline{K_m}$ and $K_{1,n}$. In Section 3, we consider some specific graphs and study the unimodality of their total domination polynomials.

\section{$\mathcal{D}_t$-classes of some graphs}

Two graphs $G$ and $H$ are said to be total dominating equivalent or simply $\mathcal{D}_t$-equivalent, if $D_t(G,x)=D_t(H,x)$ and written $G\thicksim_t H$. It is evident that the relation $\thicksim$ of $\mathcal{D}_t$-equivalent is an equivalence relation on the family $\mathcal{G}$ of graphs, and thus $\mathcal{G}$ is partitioned into equivalence classes, called the $\mathcal{D}_t$-equivalence classes. Given $G\in \mathcal{G}$, let 
\begin{center}
$[G]=\{H\in \mathcal{G} ~: H\thicksim_t G\}.$
\end{center}

If $[G]=\{G\}$, then $G$ is said to be total dominating unique or simply $\mathcal{D}_t$-unique.

It is easy to see, if two graphs $G$ and $H$ are isomorphic, then $D_t(G,x)=D_t(H,x)$, but the reverse is not always true. We have shown all graphs of order less than or equal six that are not isomorphic but have the same total domination polynomial, as Figures \ref{figeq4}, \ref{figeq5}, and \ref{f6} in Appendix. Note that all graphs of order one, two and three are $\mathcal{D}_t$-unique.   
We need the following theorems to obtain more results on $\mathcal{D}_t$-equivalence classes of some graphs: 

\begin{theorem}{\rm\cite{Dod}}
  Let $G = (V,E)$ be a graph. Then
\begin{equation*}
D_t(G, x) = D_t(G \setminus v, x) + D_t(G \odot v, x) - D_t(G \circledcirc v, x)
\end{equation*}\label{4}
where $G \odot v$ denotes the graph obtained from $G$ by removing all edges between vertices of $N(v)$ and $G \circledcirc v$ denotes the graph $G\odot v\setminus v$.
\end{theorem}

\begin{theorem}{\rm \cite{Dod}} \label{recurrence}
If  $G=(V,E)$ is  a graph and   $e=\{u,v\}\in E $ with  $N[v]=N[u]$, then 
             $D_t(G, x)=D_t(G\setminus e, x)+x^2D_t(G\setminus N[u],x)$.           
 \end{theorem}

\begin{theorem}{\rm\cite{nasrinirre}}
Let $G$ be a graph and $e=\{u,v\}$ is an edge of $G$. If $u$ and $v$ are adjacent to the support vertices, then $e$ is an irrelevant edge. That means $D_t(G,x)=D_t(G\setminus e,x)$.
\label{6}
\end{theorem}
\begin{theorem}
Let $G$ be a graph and $e=\{u,v\}\in E(G)$ is an edge of $G$ that $N_G[u]=N_G[v]$. If there is a vertex $w$ such as $N_G(w)\subseteq N_G(u)$, then $e$ is an irrelevant edge. That means $D_t(G,x)=D_t(G\setminus e,x)$.
\end{theorem}
\proof
Let $G$ be a graph and $e=\{u,v\}\in E(G)$ that $N_G[u]=N_G[v]$. By Theorem \ref{recurrence} and recording to $G\setminus N_G[u]$ has at least a isolate vertex, $w$, so $D_t(G\setminus N[u])=0$ and we have result. \qed

\begin{figure}[ht!]
 \centering
 \includegraphics[width=5.5cm, height=4.5cm]{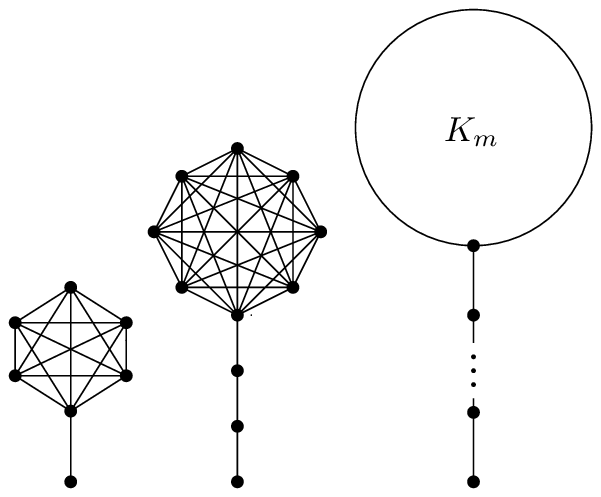}
 \caption{The lollipop graphs $L(6,1)$, $L(8,3)$ and $L(m,n)$ .}\label{figlolli}
\end{figure}

The $(m,n)$-lollipop graph is a special type of graph consisting of a complete graph $K_m$ of order $m$ and a path graph on $n$ vertices, connected with a bridge. See Figure \ref{figlolli}

\begin{corollary}\label{th21}
For every natural number $m$, the total domination polynomial of $(m,1)$-lollipop graph is equal to
\[ D_t(L(m,1),x)=x(x+1)^m-x.\]
\end{corollary}

\proof
By Theorem \ref{6}, all edges of complete graph $K_m$ in $(m,1)$-lollipop graph are irrelevant. So the total domination polynomial of this graph is equal to the total domination polynomial of the star graph $K_{1,m}$ and we have result.\qed

Generally, the total domination polynomial of $(m,n)$-lollipop graphs is obtained from the following recursive relation:
\begin{equation*}
 D_t(L(m,n),x)=xD_t(L(m,n-1),x)+x^2[D_t(L(m,n-3),x)+D_t(L(m,n-4),x)],
 \end{equation*}
 where
\[D_t(L(m,1),x)=x(x+1)^m-x,\] 
\[D_t(L(m,2),x)=x^2(x+1)^{m-1}(x+2)-(m-1)x^3-x^2,\]
\[D_t(L(m,3),x)=x^2(x+1)^{m}(x+2)-(m-1)x^4-2mx^3-2x^2,\] \[D_t(L(m,4),x)=x^2(x+1)^m(x^2+3x+1)-(m-1)x^5-2mx^4-(m+2)x^3-x^2.\]

\begin{theorem}
Let $G$ be a graph of order $n$ and  $deg v=n-1$. Then $G$ is $\mathcal{D}_t$-unique if and only if $G\setminus v$ is $\mathcal{D}_t$-unique.
\label{5}
\end{theorem}
 \proof By Theorem \ref{4}, we have
 \begin{equation*}
D_t(G, x) = D_t(G \setminus v, x) + D_t(G \odot v, x) - D_t(G \circledcirc v, x)
\end{equation*}
where $ D_t(G \odot v, x)=D_t(K_{1,n-1},x)$ and $D_t(G \circledcirc v, x)=0$. So we have result.\qed

The friendship (or Dutch-Windmill) graph $F_n$ is a graph that can be constructed by coalescence $n$
copies of the cycle graph $C_3$ of length $3$ with a common vertex. The Friendship theorem of Paul Erd\"{o}s,
Alfred R\'{e}nyi and Vera T. S\'{o}s \cite{erdos}, states that graphs with the property that every two vertices have
exactly one neighbour in common are exactly the friendship graphs.
Figure \ref{Dutch} shows some examples of friendship graphs.

\begin{figure}[ht!]
	\begin{center}
		\includegraphics[width=10cm,height=2.3cm]{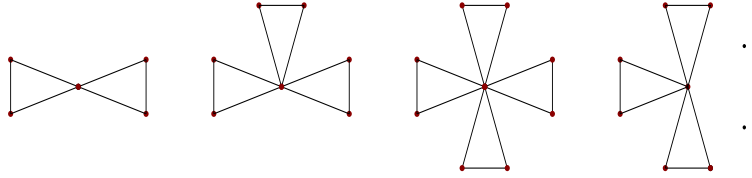}
		\caption{Friendship graphs $F_2, F_3, F_4$ and $F_n$, respectively.}
		\label{Dutch}
	\end{center}
\end{figure}

\begin{corollary}
\begin{enumerate}
\item[i)] For every $n>0$, $K_n$ is $\mathcal{D}_t$-unique.
\item[ii)] The friendship graph $F_n$ is $\mathcal{D}_t$-unique, for every $n\geq 3$.
\end{enumerate}
\end{corollary}
\proof
\begin{enumerate}
\item[i)] The result is obtained by induction and Theorem \ref{5}.
\item[ii)] By Theorem \ref{5}, since $F_n\setminus v$ is $D_t$-unique where $v$ is the center vertex of $F_n$, so we have the result.\qed
\end{enumerate}

\begin{theorem}
For every natural number $n>2$, $K_{1,n}$ is not $\mathcal{D}_t$-unique and especially
 $[K_{1,n}]\supseteq \{K_{1,n}, L(n,1), L(n,1)-e,\ldots\}$ where $e$ is any  edge of complete graph $K_n$ in lollipop graph that is not adjacent to the pendent edge of this graph. 
\end{theorem}
\proof
Let $v$ be center vertex in $K_{1,n}$. We have $D_t(K_{1,n}\setminus v,x)=0$, so $K_{1,n}\setminus v$ is not $\mathcal{D}_t$-unique and by Theorem \ref{5} we have the result. َAlso by Theorem \ref{th21} the second result is achieved.\qed

Now, we introduce an infinite family of graphs such that are total dominating equivalent with $G\circ \overline{K_m}$.
Let $G$ be a graph with vertex set $\{v_1,\ldots , v_n\}$. By $G(v_1^{m_1}, v_2^{m_2},\ldots, v_n^{m_n})$, we mean the graph obtained from $G$ by joining $m_i$ new vertices to each $v_i$, for $i = 1, \ldots , n$, where $m_1, \ldots ,m_n$ are positive integers; this graph is called sunlike. We note that by the new notation, $G\circ K_1$ is equal to $G(v_1^1, v_2^1,\ldots, v_n^1)$.

\begin{theorem}
Let $G$ be a connected graph  of order $n$. Any graphs of the family 
\[\{G\circ \overline{K_m}, (G\circ \overline{K_m})\circ \overline{K_m},  ((G\circ \overline{K_m})\circ \overline{K_m})\circ \overline{K_m},\ldots \}\] 
is not $\mathcal{D}_t$-unique.
\end{theorem}

\proof
Actually for every connected graph $G$ of order $n$, 
$$[G\circ \overline{K_m}]\supseteq \{G\circ \overline{K_m}, G(v_1^{m_1}, v_2^{m_2},\ldots, v_n^{m_n})\},$$
 where  $\sum \limits_{i=1}^{n}m_i=mn$ and for every $i$,  $m_i\geq 1$.\qed

\section{Unimodality of total domination polynomial}


In this section, we consider some specific graphs and study the unimodality of their total domination polynomials. We think that the total domination polynomial of a graph is unimodal \cite{nasrin}. 
We need the following theorem to state and prove some new results for the unimodality of the total domination polynomial of graphs. 

\begin{theorem}\label{log} {\rm\cite{log}} 
Let $f(x)$ and $g(x)$ be polynomials with positive coefficients.
\begin{enumerate}
\item[i)]If both $f(x)$ and $g(x)$ are log-concave, then so is their product $f(x)g(x)$.
\item[ii)]If $f(x)$ is log-concave and $g(x)$ is unimodal, then their product $f(x)g(x)$ is unimodal.
\item[iii)]If both $f(x)$ and $g(x)$ are symmetric and unimodal, then so is their product $f(x)g(x)$.
\end{enumerate}
\end{theorem}
 
So if polynomials $P_i(x)$ for $i = 1,2, \ldots, n$ with positive coefficients are log-concave, then $\prod_{k=1}^n P_k(x)$ is log-concave as well.
Here, we introduce a family of graphs whose total domination polynomial are unimodal.

\begin{figure}
\begin{center}
\includegraphics[width=9.75cm,height=2.6cm]{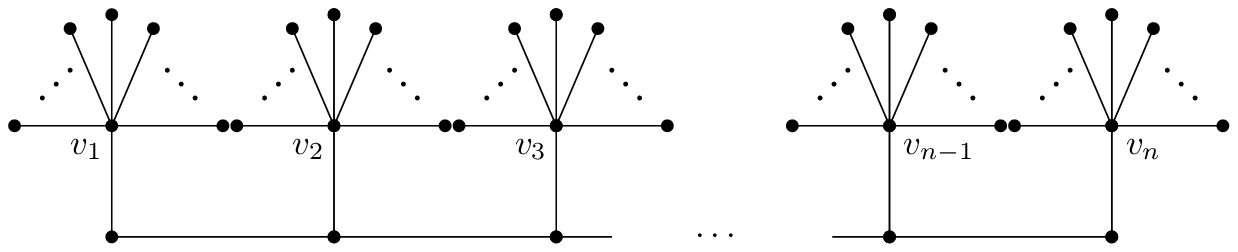}
\includegraphics[width=9.25cm, height=2.25cm]{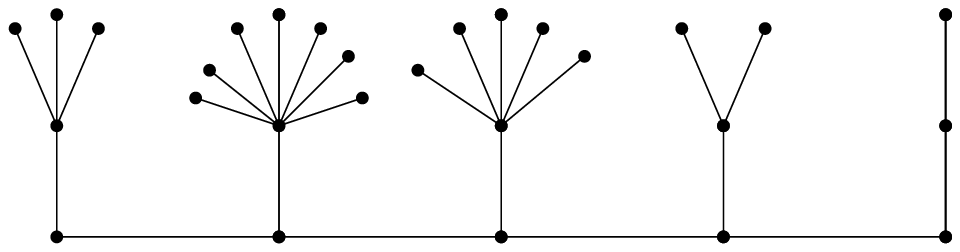}
\caption{The graph $F(n,k)$ and $F(5,9,7,4,3)$.}
\label{FNK}
\end{center}
\end{figure}

An $(n,k)$-firecracker $F(n,k)$ is a graph obtained by the concatenation of $n$, $k$-stars $S_k$ by linking one leaf from each. Also we generalize the definition of firecracker graphs. An $(k_1,k_2,\ldots,k_n)$-firecracker $F(k_1,\ldots, k_n)$ is a a graph obtained by the concatenation of $k_i$-stars $S_{k_i}$ by linking one leaf from each (see Figure \ref{FNK}). 
  
\begin{theorem}{\rm \cite{nasrinirre}}
	\begin{enumerate} 
		\item[(i)] 
		For every natural numbers $n$ and $k\geq 3$,  $D_t(F(n,k),x)=(x(x+1)^{(k-1)}-x)^n$.
 \item[(ii)]  $D_t(F(k_1,\ldots, k_n))=\prod \limits_{i=1}^n (x(x+1)^{(k_i-1)}-x)$. 
 \end{enumerate} 
\end{theorem}

\begin{theorem}
For natural numbers $n$, $k\geq 3$ and $k_i\geq 3$ ($1\leq i\leq n$) the total domination polynomial of graphs $L(n,1)$, $F(n,k)$ and $F(k_1,k_2,\ldots,k_n)$ are  unimodal.
\end{theorem}

\proof
Total domination polynomial of these graphs is equal to product of the total domination polynomial of some star graphs.
For every natural number $n$ , $D_t(K_{1,n},x)=x(x+1)^n-x$ is unimodal and in particular log-concave. So by Theorem \ref{log} we have results.\qed

The following theorem gives us many graphs whose total domination polynomials are unimodal:

\begin{theorem}
Let $G$ be a graph of order $n$ with $r$ isolated vertices. Total domination polynomial of every graph of family
 \[\{G\circ \overline{K_m}, (G\circ \overline{K_m})\circ \overline{K_m}, ((G\circ \overline{K_m})\circ \overline{K_m})\circ \overline{K_m},\cdots \}\] 
 is unimodal.
\end{theorem}

\proof
For every graph of order $n$ with $r$ isolated vertices we have
\[D_t(G\circ \overline{K_m},x)=x^n(1+x)^{m(n-r)}[(x+1)^n-1]^r.\]
So by Theorem \ref{log} the total domination polynomial of these graphs is log-concave and so  is unimodal.\qed

The generalized friendship graph $F_{n,q}$ is a collection of $n$ cycles (all of order $q$), meeting at a common vertex (see Figure \ref{Dutch}). The generalized friendship graph may also be referred to as a flower (\cite{Schiermeyer}). For $q=3$ the graph $F_{n,q}$ is denoted simply by $F_n$ and is friendship graph as known.

\begin{figure}[ht]
\begin{center}
\includegraphics[width=5.5in]{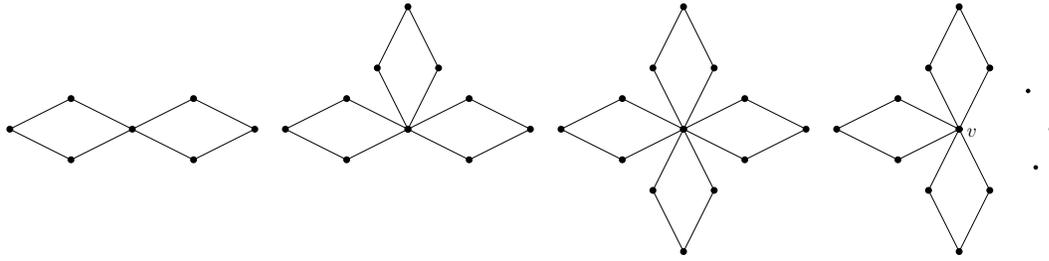}
\end{center}
\vspace{-.5cm}
\caption{Generalized friendship graphs $F_{2,4}$, $F_{3,4}$, $F_{4,4}$ and $F_{n,4}$.}\label{Dutch}
\end{figure}

\begin{figure}[!ht]
	\hspace{4cm}
	\includegraphics[width=6.5cm,height=2cm]{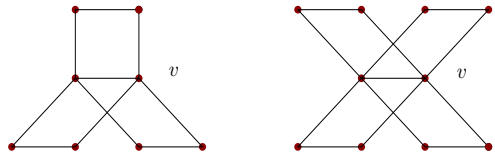}
	\caption{ \label{figure6} The book graphs $B_3$ and $B_4$, respectively.}
\end{figure}

The $n$-book graph $B_n$ can be constructed by bonding $n$ copies of the cycle graph $C_4$ along a common edge $\{u, v\}$, see Figure \ref{figure6}. Here we compute the total domination polynomial of   $n$-book graphs. 

\begin{theorem}
For each natural number $n$,
\begin{enumerate}
\item[i)]Total domination polynomial of $n$-book graph $B_n$ is unimodal.
\item[ii)]Total domination polynomial of graph $F_{n,4}$ is unimodal.
\end{enumerate}
\end{theorem}

\proof
\begin{enumerate}
\item[i)] For each natural number $n$ we have $D_t(B_n,x)=(x(x+1)^n+x^n)^2$ (see \cite{nasrin})
\[D_t(B_n,x)=(x^{n+1}+(n+1)x^n+\binom{n}{n-2}x^{n-1}\ldots +\binom{n}{2}x^3+nx^2)^2.\]
By Theorem \ref{log} this polynomial is log-concave and so unimodal.
\item[ii)] Since  $D_t(F_{n,4},x)=x^{n+1}(x+2)^n[(x+1)^n+x^{n-1}]$ (see \cite{nasrin3}) and 
\[\big(2^i\binom{n}{i}\big)^2\geq 2^{i-1}\binom{n}{i-1}2^{i+1}\binom{n}{i+1}=2^{2i}\binom{n}{i-1}\binom{n}{i+1},\]
so this polynomial is unimodal.\qed
\end{enumerate}

Some results about unimodality of polynomials prove by position of their roots.
\begin{theorem}\label{prealp}{\rm \cite{Comtet}}
If a polynomial $p(x)$ with positive coefficients has all real roots, then is log-concave and unimodal.
\end{theorem}

Here, we introduce some family of graphs whose roots of total domination polynomial are real, So their total domination polynomial are log-concave and unimodal.
\\
The helm graph $H_n$ is obtained from the wheel graph $W_n$ by attaching a pendent edge at each vertex of the $n$-cycle of the wheel. We define generalized helm graph $H_{n,m}$, the graph is obtained from the wheel graph $W_n$ by attaching $m$ pendent edges at each vertex of the $n$-cycle of the wheel.
We recall that corona product of two graphs $G$ and $H$ is denoted by $G\circ H$  and was introduced by Harary \cite{Harary, Harar}. This graph formed from one copy of $G$ and $|V(G)|$ copies of $H$, where the $i$-th vertex of $G$ is adjacent to every vertex in the $i$-th copy of $H$. 

\begin{figure}[ht]
 \centering
 \includegraphics[height=3.4cm]{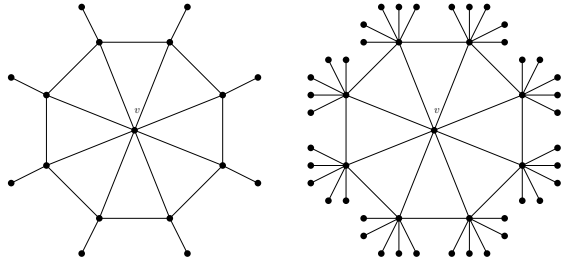}
 \caption{Helm graph $H_8$  and generalized helm graph $H_{8,5}$, respectively.}\label{fighelm}
\end{figure}

\begin{theorem}{\rm \cite{nasrin3}}
For natural numbers $m,n$, $D_t(H_{n,m},x)=x^n(x+1)^{mn+1}$, specially for $m=1$ we have $D_t(H_{n},x)=x^n(x+1)^{n+1}$
\end{theorem}

\begin{figure}
	\begin{center}
		\includegraphics[width=3.7cm,height=3cm]{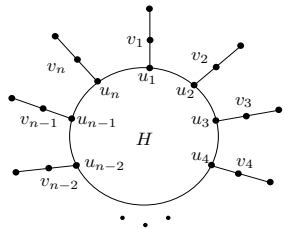}
		\caption{The graph $H(3)$.}
		\label{H(3)}
	\end{center}
\end{figure}

By the definition, the graph $H(3)$ is obtained by identifying each vertex of $H$ with an end vertex of a $P_3$ (\cite{Nasrin2}). See Figure \ref{H(3)}.

\begin{theorem}\label{H}{\rm \cite{nasrinirre}}
	For any graph $H$ of order $n$, we have $D_t(H(3),x)=x^{2n}(x+2)^n$.
\end{theorem}

 By Theorem \ref{prealp} the total domination polynomial of all graphs with real roots such as $H_n$, $H_{m,n}$, $H(3)$ and sunlike graphs $G(v_1^{k_1},v_2^{k_2},\ldots, v_n^{k_n})$  are  unimodal.

\newpage
\section{Appendix: $\mathcal{D}_t$-equivalence classes of graphs of order $\leq 6$} 
\begin{figure}[ht!]
	\centering
	\includegraphics[width=5cm, height=3.5cm]{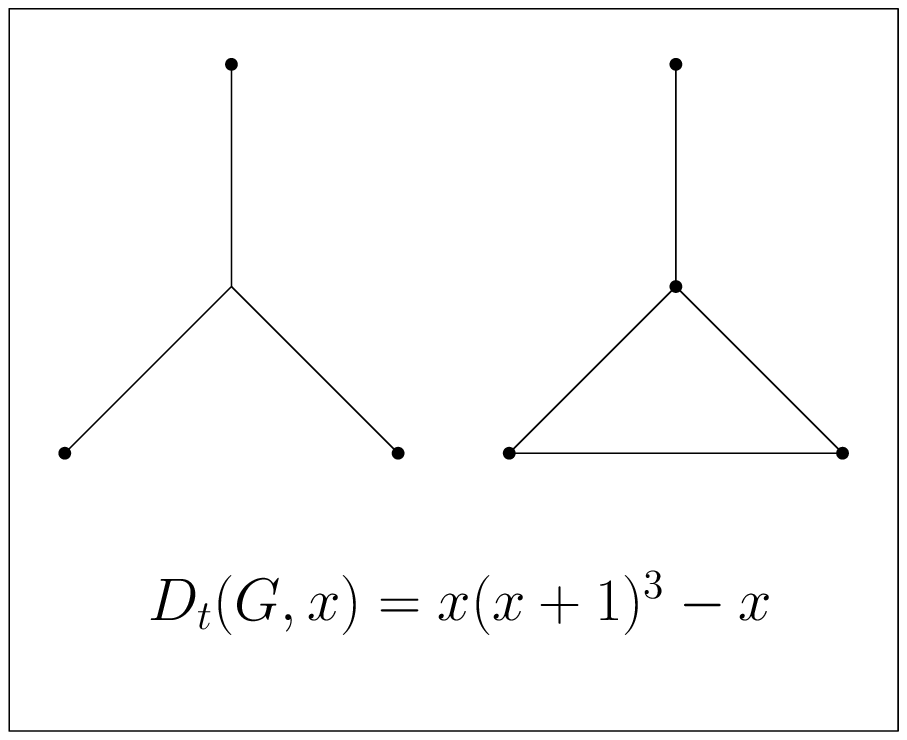}
	\caption{ The $\mathcal{D}_t$-equivalence class of connected graphs of order $4$. }\label{figeq4}
\end{figure}
\begin{figure}[ht!]
	\centering
	\includegraphics[width=11.5cm, height=9cm]{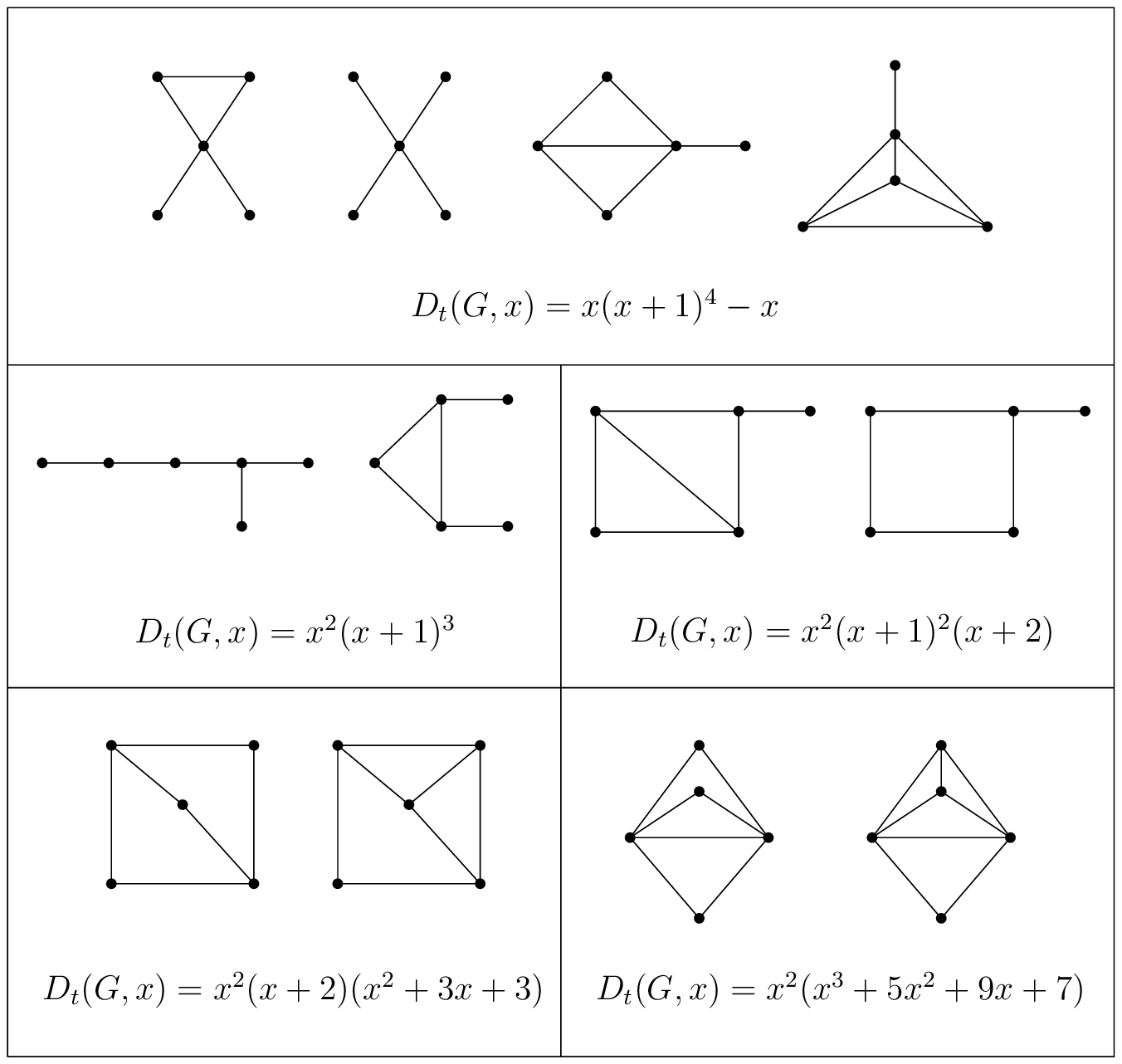}
	\caption{ The $\mathcal{D}_t$-equivalence classes of connected graphs of order $5$. }\label{figeq5}
\end{figure}

\begin{figure}[ht!]
	\centering
	\includegraphics[width=14cm, height=18cm]{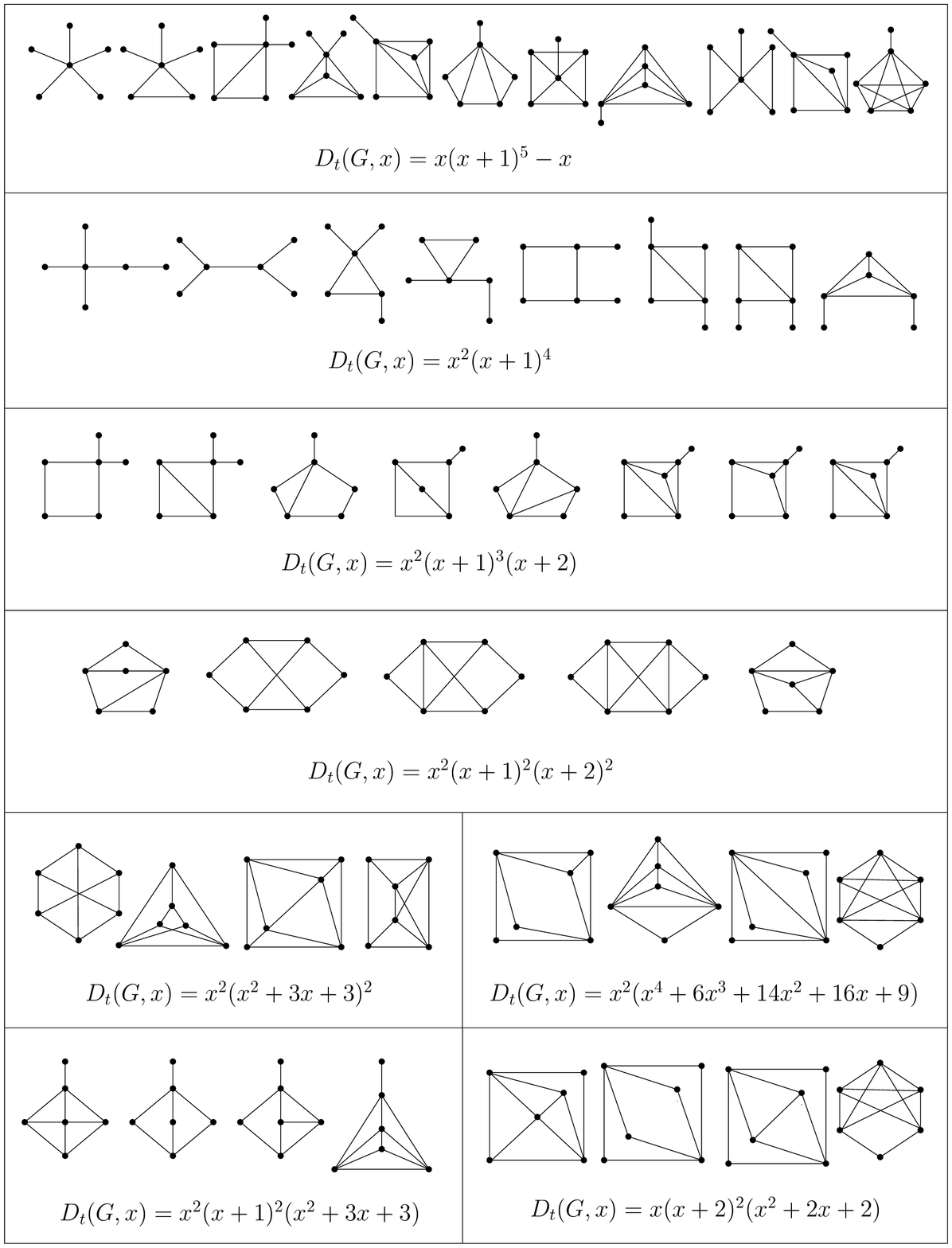}
\end{figure}
\begin{figure}[ht!]
	\centering
	\includegraphics[width=14cm, height=18cm]{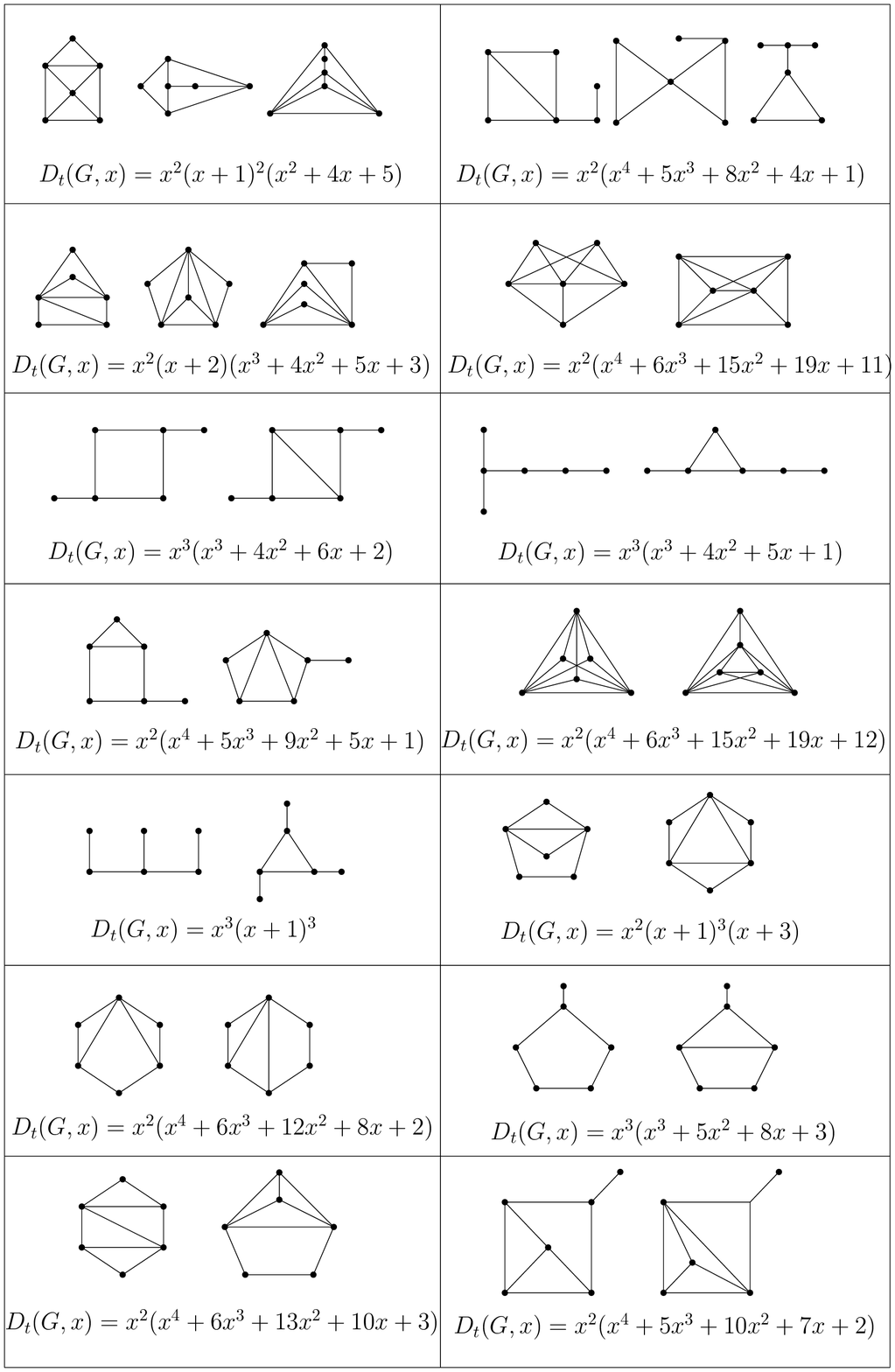}
	\caption{ The $\mathcal{D}_t$-equivalence classes of connected graphs of order $6$.}\label{f6}
\end{figure}

\end{document}